\documentclass[letterpaper, 10 pt, conference]{ieeeconf}  % Comment this line out if you need a4paper

\IEEEoverridecommandlockouts                              % This command is only needed if 
\let\labelindent\relax

\usepackage{amsmath,amssymb,amsfonts,amsthm}
\usepackage[scr=dutchcal]{mathalpha}
\usepackage{grffile}
\usepackage{multicol}
\usepackage{caption}
\usepackage{url}
\usepackage{enumitem}

\newcommand{\mcl}[1]{\mathcal{ #1}}
\newcommand{\mbf}[1]{\mathbf{ #1}}
\newcommand{\norm}[1]{\left\Vert #1\right\Vert}
\newcommand{\hinf}{\ensuremath{H_{\infty}}}
\newcommand{\ip}[2]{\left\langle{#1},{#2}\right\rangle}
\newcommand{\half}{\frac{1}{2}}

\newcommand{\bmat}[1]{\begin{bmatrix} #1\end{bmatrix}}

\newcommand{\mat}[1]{\begin{matrix}#1\end{matrix}}

\newcommand{\R}{\mathbb{R}}

\newcommand{\Z}{\mathbb{Z}}

\newtheorem{thm}{Theorem}
\newtheorem{defn}[thm]{Definition}

\newtheorem{cor}[thm]{Corollary}

\newtheorem{ex}[thm]{\textbf{Example}}

\newcommand{\PI}{\pmb{\Pi}}
\newcommand{\pie}{\scalebox{0.9}[1.2]{$\mathit{\Pi}$}}

\newcommand{\fourpi}[4]{\hspace{0.5mm}\pie\hspace{-1.mm}\left[\footnotesize\begin{array}{c|c}
#1&#2\\\hline #3 & \{#4\}
\end{array}\right]}

\allowdisplaybreaks

\title{\LARGE \bf
A Computational Method for $H_2$-optimal Estimator and State Feedback Controller Synthesis for PDEs
}

\author{Sachin Shivakumar$^{1}$ and Matthew M. Peet$^{2}$% <-this % stops a space
\thanks{This work was supported by the National Science Foundation under grants No. 1739990 and 1935453}% <-this % stops a space
\thanks{$^{1}$ Sachin Shivakumar\{{\tt\small sshivak8@asu.edu}\} and Matthew M. Peet\{{\tt\small mpeet@asu.edu}\} are with School for Engineering of Matter, Transport and Energy, Arizona State University, USA
        }%
}

\begin{document}

\maketitle
\thispagestyle{empty}
\pagestyle{empty}

%%%%%%%%%%%%%%%%%%%%%%%%%%%%%%%%%%%%%%%%%%%%%%%%%%%%%%%%%%%%%%%%%%%%%%%%%%%%%%%%
\begin{abstract}
In this paper, we present solvable, convex formulations of $H_2$-optimal state estimation and state-feedback control problems for a general class of linear Partial Differential Equations (PDEs) with one spatial dimension. These convex formulations are derived by using an analysis and control framework called the `Partial Integral Equation' (PIE) framework, which utilizes the PIE representation of infinite-dimensional systems. Since PIEs are parameterized by Partial Integral (PI) operators that form an algebra, $H_2$-optimal estimation and control problems for PIEs can be formulated as Linear PI Inequalities (LPIs). Furthermore, if a PDE admits a PIE representation, then the stability and $H_2$ performance of the PIE system implies that of the PDE system. Consequently, the optimal estimator and controller obtained for a PIE using LPIs provide the same stability and performance when applied to the corresponding PDE. These LPI optimization problems can be solved computationally using semi-definite programming solvers because such problems can be formulated using Linear Matrix Inequalities by using positive matrices to parameterize a cone of positive PI operators. We illustrate the application of these methods by constructing observers and controllers for some standard PDE examples.
\end{abstract}

\section{Introduction}

Various metrics are used to design observers and controllers for dynamical systems such that the performance with respect to the metric is optimal. For example, some standard choice of metrics are quadratic cost function on state and inputs, $H_{\infty}$-norm, $H_2$-norm, etc. Among these properties, this paper focuses on optimal estimation and control problems with $H_2$-norm of the system as the metric because $H_2$-norm gives information on the system's output behavior in the presence of an impulse input, white noise input, or initial conditions. Furthermore, $H_2$-optimal control problems are also a generalization of standard control problems, such as Linear Quadratic Regulator (LQR) and Linear Quadratic Gaussian (LQG).

The gains of the observer and controller for the closed-loop Partial Differential Equation (PDE) that provide optimal $H_2$-performance have important applications in systems that experience impulsive and stochastic noise inputs. Alternatively, such optimal estimation/control problems also appear as sub-problems in applications such as sensor/actuator placement to improve closed-loop performance. For example, \cite{armaou2006optimal} deals with determining optimal sensor/actuator location along with the observer and controller gains by solving a nonlinear optimization problem (Also, see \cite{munz2014sensor} for an iterative convex optimization formulation). Such simultaneous optimization problems (finding sensor/actuator location and observer/controller gains) have also appeared in applications such a beam vibration control \cite{liu2006computational,ambrosio2012h2}, flow control \cite{chen2011h2}, etc.

Despite the importance and various applications, finding the estimator/controller that optimizes the $H_2$-norm of the closed-loop PDE system is difficult. This is because there are no methods that are non-conservative and can find $H_2$-optimal estimators and controllers with provable properties. This is primarily due to the inability to obtain provable $H_2$-norm bounds for most PDE systems, barring some specific cases where an analytical expression for the transfer function can be obtained \cite{tf_PDE} or when the PDE has finite number of unstable modes \cite{bamieh2002distributed,lasiecka2000control}. Since such PDEs are rare, most of the existing methods rely on early or late-lumping approximation of the PDE to find $H_2$-optimal observers and controllers --- providing observer/controller gains that have provable performance only when approximation errors are zero.

Early-lumping approaches, such as \cite{mechhoud2014estimation}, discretize the PDE solution space (using Galerkin projection, spatial discretization or modal decomposition), to obtain an Ordinary Differential Equation (ODE) approximation of the PDE. These approaches, however, do not have provable bounds on $H_2$-norm because approximating a PDE by an ODE lead to approximation errors, as a result of which there may be little or no relationship between the solutions of the original problem and its approximation. In short, the bounds on $H_2$-norm are \textit{unprovable} in the sense that the bounds so obtained can exceed or fall short of the true $H_2$-norm while not providing any metric to determine the accuracy of the bound.

On the other hand, late-lumping approaches, such as \cite{schmidt2016reduced,benner2018numerical,moghadam2013boundary} (for Time delay systems, see \cite{zhang2006optimal}), do not approximate the PDE, however, require solving an Operator-valued Ricatti equation (ORE) with infinitedimensional, possibly unbounded, operators -- e.g., ORE for LQR control of the heat equation woudl have second-order spatial differential operators. Since the set of such operators is not algebraically closed, one cannot easily parameterize the unknown operators to solve these OREs. Thus, the operators are projected onto a finite-dimensional subspace prior to solving, which often lead to conservative bounds on the $H_2$-norm (or \textit{unprovable} bounds, if the truncation errors from projection cannot be bounded).

The goal of this paper is to solve the problem of `finding $H_2$-optimal estimator/controller for PDEs' while overcoming the shortcomings of the existing methods, namely, conservatism or \textit{unprovability} of the bounds on $H_2$-norm. A method was proposed in \cite{danilo} to provably upper-bound the $H_2$-norm of a PDE system by utilizing an alternative representation for PDEs called the Partial Integral Equation (PIE) representation --- an equivalent representation defined by bounded linear operators called Partial Integral (PI) operators. 

The primary motivation behind the PIE representation is the lack of a universal parameterization or representation of analysis and control problems for PDE systems -- an inconvenience that is further compounded by the presence of unbounded operators and auxiliary constraints. Although \cite{shivakumar_representation_TAC} introduces a general parametric form for a large class of PDEs, the number of parameters changes with the order of derivatives in the PDE, the boundary conditions, etc. To overcome this issue of varying number of parameters, the PIE representation is used that has a fixed number of parameters given by at most $12$ Partial Integral (PI) operators. Subsequently, it was proved that the analysis and control problems for PIEs can be formulated are solvable convex optimization problems with PI operator-valued constraints called Linear PI Inequalities (LPIs). For e.g., exponential stability and input-to-output $L_2$-gain \cite{shivakumar_representation_TAC}, $H_{\infty}$-optimal controller synthesis \cite{shivakumar2020duality}, etc. 

Inspired by the success of the PIE framework in other analyses and control problems of PDEs, \cite{danilo} used the PIE representation and the LPIs to obtain upper-bounds on the $H_2$-norm for PDE and Time-Delay Systems (TDSs). Although the method proposed therein was not conservative and typically outperformed other existing methods, the results were limited to finding the $H_2$-norm and a generalized observability Gramian. While the results from \cite{danilo} can indeed be extended to construct state-estimators with $H_2$-optimal performance of the error system, one cannot use the results for controller synthesis. Since the LPI optimization approach uses parametrized Lyapunov functions to search for observer and controller with provable properties for PIE systems, one often runs into the bilinearity issue where the constraints of the optimization problem are bilinear in the Lyapunov function parameters and gains of the observer (or controller). In the case of $H_2$-optimal state estimation problem, one can overcome the bilinearity via an invertible variable change. However, the same technique cannot be used in the control problem because the such an invertible variable change does not exist --- analogous to the ODE case. Thus, to resolve the bilinearity in the $H_2$-optimal control problem, one must employ the dual PIE representation of a PDE as was done in \cite{shivakumar2020duality}. Although \cite{shivakumar2020duality} proves that the $\hinf$-norm of a PIE and its dual are equal, there are no results on $H_2$-norm equivalence. 

Using input-to-output charactertization of $H_2$-norm one might claim the equivalence of $H_2$-norm, however, such a charactertization depends on transfer functions --- a notion that relies on the existence of well-defined solution map of the system. In the case of PDEs, proving existence such a map and finding the transfer function is often challenging. Thus, to avoid this reliance of proving existence and uniqueness, we use a state-to-output characterization of $H_2$-norm that is identical to the input-to-output characterization when the solution map of a PDE is a strongly continous semigroup (see \cite{danilo} for details). Hence, to utilize the dual PIE representation for controller synthesis, we will first establish equivalence of $H_2$-norm of a PIE and its associated dual in Sec.~\ref{sec:LPI}. Once we establish this equivalence, we can pose the controller synthesis problem for the dual PIE (instead of the original PIE) and convexify the bilinear constraints via an invertible variable change --- similar to the approach used in $\hinf$-optimal state-feedback control in \cite{shivakumar2020duality}. 

To summarize, we will use PIE representation and dual PIE representation presented in Secs.~\ref{subsec:PIE} and \ref{subsec:h2norm}, respectively, to propose convex optimization formulations of $H_2$-optimal state estimation and state-feedback control problems for a general class of linear PDEs in one spatial dimension presented in Sec.~\ref{subsec:PDE}. These formulations are presented in the form of LPI optimization problems in Sec.~\ref{sec:LPI}, which can be solved using semi-definite programming or PIETOOLS \cite{shivakumar_2020ACC}, an open source MATLAB toolbox. We apply the results of Sec.~\ref{sec:LPI} to numerical examples in Sec.~\ref{sec:numerical}.

\section{Preliminaries}\label{sec:preliminaries}
In this section, we will introduce the notation used in the paper along with some rudimentary information on the PIE framework used in the analysis and control of PDEs. We also describe, through an example, the type of PDEs for which the framework is applicable. 
\subsection{Notation}
We denote the set of all real Lesbegue square-integrable functions on the spatial domain $[a,b]\subset \R$ as $L_2^n[a,b]$. Similarly $L_2^n[0,\infty)$ denotes square-integrable real-valued signals where $[0,\infty)$ is a temporal domain.  
% $W_{k}^n[a,b]$ denotes the Sobolev space \begin{align*}
% W_{k}^n[a,b] := \{\mbf x \in L_2^n[a,b] : \partial^i_s \mbf x\in L_2^n[a,b] ~\text{for all} ~i\le k\}
% \end{align*}
% where $\partial_s^j \mbf x$ denotes the partial derivative $\frac{\partial^j\mbf x}{\partial s^j}$. $W_{k}^n[a,b]$ is equipped with the standard Sobolev inner-product, denoted $\ip{\cdot}{\cdot}_{W_k}$. In this paper, Sobolev spaces are associated only with a compact spatial domain.
For brevity, $\R L_2^{m,n}[a,b]$ denotes the space $\R^{m}\times L_2^{n}[a,b]$ where for $\bmat{x_1\\\mbf x_2}, \bmat{y_1\\\mbf y_2}\in \R L_2^{m,n}$, the inner-product is defined as
\begin{align*}
\ip{\bmat{x_1\\\mbf x_2}}{\bmat{y_1\\\mbf y_2}}_{\R L_2} := x_1^Ty_1 + \ip{\mbf x_2}{\mbf y_2}_{L_2}.
\end{align*}
Occasionally, we omit the domain and simply write $L_2^n$ or $\R L_2^{m,n}$. We also omit the inner-product subscripts whenever it is clear from the context. For functions of time and space, e.g. $\mbf x(t,s)$,  $\dot{\mbf x}$ is used to denote the partial derivative $\frac{\partial \mbf x}{\partial t}$ and $\partial_s \mbf x$ to denote $\frac{\partial \mbf x}{\partial s}$.

We use the bold font to indicate functions of space (or time and space), e.g. $\mbf x \in L_2^n[a,b]$ and calligraphic font, e.g. $\mcl{A}$, to represent bounded linear operators on Hilbert spaces, e.g. $\mcl A \in B(X)$ where $X$ is a Hilbert space with inner product $\ip{\cdot}{\cdot}_X$. For any $\mcl A\in B(X)$, $\mcl A^*$ denotes the adjoint operator satisfying $\ip{\mbf x}{\mcl A\mbf y}_{X} = \ip{\mcl A^*\mbf x}{\mbf y}_{X}$ for all $\mbf x, \mbf y \in X$. Lastly, we define the truncation operator $P_T$ as
\[
P_T z(t) = \left\lbrace \mat{z(t), &t\le T\\ 0, &\text{otherwise,}}\right.\quad z\in L_{2e}[\R^+].
\]

\subsection{Partial Integral Equations}\label{subsec:PIE}
Having introduced some standard notation, we will next introduce some notation for PI operators that are used to define PIE systems. PI operators are bounded linear operators on $\R L_2$ and are elements of the $^*$-algebra $\PI_4$, which is defined as follows.
\begin{defn}\label{def:4PI}
We say $\mcl P \in \PI_4 \subset B(\R L_2^{m_1,n_1},\R L_2^{m_2,n_2})$ if there exists a matrix $P$ and polynomials $Q_1,Q_2,R_0,R_1$, and $R_2$ such that
\begin{align*}
&\mcl P=\fourpi{P}{Q_1}{Q_2}{R_i}\bmat{x\\\mbf{y}}(s) := \bmat{Px + \int_{a}^{b}Q_1(s)\mbf{y}(s)ds\\Q_2(s)x+ \mcl R\mbf{y} (s)},\\
&\left(\mcl R\mbf y\right)(s)\hspace{-1.5mm}= \hspace{-1.5mm}R_0(s) \mbf y(s) +\hspace{-1.5mm}\int\limits_{a}^s  \hspace{-1.5mm}R_1(s,\theta)\mbf y(\theta)d \theta+\hspace{-1.5mm}\int\limits_s^b \hspace{-1.5mm}R_2(s,\theta)\mbf y(\theta)d \theta.
\end{align*}
\end{defn}
Since the set of PI operators is closed under composition, addition, and adjoint, explicit formulae for these operations on PI operators can be obtained in terms of operations on the parameters; See  \cite{shivakumar_representation_TAC} for details and formulae.

The notation $\fourpi{P}{Q_1}{Q_2}{R_i}$ is used to indicate the PI operator associated with the matrix $P$ and polynomial parameters $Q_i$, $R_j$. The dimensions ($m_1,n_1,m_2,n_2$) are inherited from the dimensions of the matrices $P\in \R^{n_2 \times n_1}$ and polynomials $R_0(s)\in \R^{m_2 \times m_1}$. When clear from context, we will omit the dimensions of the domain and range and simply use $\R L_2$. In the case where a dimension is zero, we use $\emptyset$ in place of the associated parameter with zero dimension. For example, if $m_1=0$, we have an operator of the form
\[
\fourpi{\emptyset}{\emptyset}{Q_2}{R_i}.
\]

In this paper, we consider estimation and control problems for the PIE systems of the form 
\begin{align}\label{eq:PIE_full}
\bmat{\mcl{T}\dot{\mbf{x}}(t)\\z(t)}  = \bmat{\mcl{A}&\mcl{B}\\\mcl C& D }\bmat{\mbf x(t)\\w(t)},
\end{align}
 where $\mcl T$ acts as a filter on the PIE state $\mbf x$, $\mcl A$ is the generator, $\mcl B$ determines the influence of input $w$ on the dynamics, $\mcl C$ determines the state contribution to the output $z$, and $D$ determines the input contribution to $z$. For systems with finite-dimensional inputs and outputs, $\mcl T, \mcl A, \mcl B, \mcl C$ are PI operators, and $D$ is a matrix. We can define a set of minimal requirements that a solution for a PIE of the above form must satisfy as shown below.
 \begin{defn}\label{def:pie}
        Given PI operators $\mcl T$, $\mcl A$, $\mcl B$, $\mcl C$ and matrix $D$, we say $\{\mbf x, z\}$ satisfy the system $\Sigma(\mcl T, \mcl A, \mcl B,\mcl C, D)$ for initial conditions $\mbf x_0\in \R L_2$ and inputs $w\in L_{2e}[\R^+]$, if $\mcl T\mbf x$ is differentiable for all $t\in \R^+$, $z\in L_{2e}[\R^+]$, and $\mbf x$ and $z$ satisfy the equations
        \begin{align}\label{eq:PIE_general}
        \bmat{\mcl{T}\dot{\mbf{x}}(t)\\z(t)}  = \bmat{\mcl{A}&\mcl{B}\\\mcl C& D }\bmat{\mbf x(t)\\w(t)},
        \end{align}
        for all $t\in \R^+$ and $(\mcl T\mbf x)(0) = \mbf x_0$.
 \end{defn}

 \noindent \textbf{Notation:} Note that in Def.~\ref{def:pie}, we use the notation $\Sigma(\mcl T, \mcl A, \mcl B,\mcl C, D)$ to denote a PIE system. Extending this notation to special case of PIEs with one or more of the null parameters, we will simply use $`-'$ in the corresponding location. For example, $\Sigma(\mcl T, \mcl A,-,\mcl C, -)$ implies $\mcl B=0$ and $D = 0$ --- i.e., there are no inputs.

 \subsection{GPDE: A generalized class of PDEs}\label{subsec:PDE}
The class of PDEs considered that may admit a PIE representation includes PDEs that have: ODE coupling, $n$th order spatial derivatives, integral terms, boundary terms, inputs, and outputs. Using the differential operator notation $\mbf D^i\mbf x:= \text{col}(\partial_s^0\mbf x,\cdots,\partial_s^i\mbf x)$, Dirac operator $\Delta_a \mbf x:= \mbf x(a)$, and integral operator $\mbf I_{[a,b]}\mbf x:=\int_a^b \mbf x(s)ds$, we can compactly represent such PDEs as
\begin{align}\label{eq:PDE-full}
    &\bmat{\dot{x}\\\dot{\mbf x}} = A_{xwu}\bmat{x\\w\\u}+\fourpi{\emptyset}{B_{x,\mbf x}}{\emptyset}{A_i}\mbf D^n\mbf x+B_{b\mbf x}\bmat{\Delta_a\\\Delta_b}\mbf D^{n-1}\mbf x,\notag\\
    &\bmat{z\\y} = C_{xwu}\bmat{x\\w\\u}+C_b\bmat{\Delta_a\\\Delta_b}\mbf D^{n-1}\mbf x+\mbf I_{[a,b]}C_{\mbf x}\mbf D^n\mbf x,\notag\\
    &B_{b}\bmat{\Delta_a\\\Delta_b}\mbf D^{n-1}\mbf x +\mbf I_{[a,b]}B_I\mbf D^n \mbf x= B_{xwu}\bmat{x\\w\\u},
\end{align}
where $x$ is an ODE state (function of time) and $\mbf x$ is $n$-times spatially differentiable PDE state (function of time and space). Following the typical $9$-matrix representation of ODEs \cite{boyd}, we separate outputs into two categories: observed $y$ and regulated $z$. Likewise, inputs are separated into disturbance $w$ and control $u$. The above parameterization is general in the following sense:
\begin{itemize}
    \item it allows $n$th-order derivatives $\mbf D^n$, boundary terms $\Delta$, and integral terms $\mbf I$, all of which can impact the dynamics (via $A_i$ and $B_{b\mbf x}$), the outputs (via $C_b$ and $C_{\mbf x}$), and boundary conditions (via $B_b$ and $B_I$).
    \item it allows coupling with ODE and influence of inputs via $B_{x,\mbf x}$, $B_{b\mbf x}$, $A_{xwu}$, and $B_{xwu}$.
\end{itemize}

In fact, this parameterization can be further extended to PDEs with mixed orders of differentiability. While we will refer to~\cite{shivakumar_representation_TAC} for the full class of linear PDEs that admit such a PIE representation Eq.~\eqref{eq:PIE_full}, we include the following example for illustration. 

\begin{ex}\label{ex:EB-representation}
Consider the vibration suppression problem for a cantilevered Euler-Bernoulli beam
\begin{align*}
&\dot{\mbf{x}}(t,s) = \bmat{0&-0.1\\1&0}\partial_s^2\mbf{x}(t,s)+\bmat{1\\0}w(t)+\bmat{1\\0}u(t),\\
&\bmat{1&0}\mbf x(t,0) = \bmat{1&0}\partial_s\mbf x(t,0) = 0,\\
&\bmat{0&1}\mbf x(t,1) =\bmat{0&1}\partial_s\mbf x(t,1) = 0,
\end{align*}
where we define the state as $\mbf{x} = col(\partial_t\eta, \partial_s^2\eta)$ where $\eta$ is displacement, $w$ is external disturbance and $u$ is control input. To regulate a combination of vibrations and control effort we defined $z(t) = \bmat{\int_0^1 \eta(t,s) ds&u(t)}^T$. The goal is to find the controller gains $\mcl K \,:\,\mbf x(t)\mapsto u(t)$ that minimizes a certain objective ($\hinf$-norm, $H_2$-norm, LQR cost function, etc.). For this PDE, we can find a PIE representation $\Sigma(\mcl T,\mcl A, \bmat{\mcl B_1\; \mcl B_2}, \mcl C, \bmat{D_1\; D_2})$ where $\mcl B_1$ and $\mcl B_2$ are the parameters corresponding to the influence of the disturbance $w$ and control $u$ on the dynamics. Likewise, $D_1$ and $D_2$ correspond to the contributions of $w$ and $u$ to the output $z$.
To illustrate, we derive the PIE representation from Cauchy's rule for repeated integration, which gives us the identity
\[
\mbf x(s)=\mbf x(0)+s\partial_s\mbf x(0)+\int_0^s (s-\theta)\partial_s^2\mbf x(\theta)d\theta.
\]
Substituting the boundary conditions, we obtain the direct relationship
\begin{align*}
&\mbf x(s) = \\
&\int_0^s \bmat{(s-\theta)&\hspace{-0.5mm}0\\0&\hspace{-0.5mm}0}\partial_s^2\mbf x(\theta)d\theta + \int_s^1 \bmat{0&\hspace{-0.5mm}0\\0&\hspace{-0.5mm}(\theta-s)}\partial_s^2\mbf x(\theta)d\theta.
\end{align*}
Substituting this expression into the dynamics and denoting $\underline{\mbf x}:=\partial_s^2\mbf x$, we obtain the PIE representation
\begin{align*}
&\partial_t\left(\mat{\int_0^s \bmat{(s-\theta)&0\\0&0}\underline{\mbf x}(t,\theta)d\theta &\\ &\hspace{-40mm}+\int_s^1 \bmat{0&0\\0&(\theta-s)}\underline{\mbf x}(t,\theta)d\theta}\right) \\
&\quad = \bmat{0&-0.1\\1&0}\underline{\mbf{x}}(t,s)+\bmat{1\\0}w(t)+\bmat{1\\0}u(t).
\end{align*}
Finally, by inspection, we identify the non-zero parameters in the Partial Integral operators $\mcl{T, A}, \mcl B_i, \mcl C, D_i$ as
\begin{align*}
&\mcl T=\fourpi{\emptyset}{\emptyset}{\emptyset}{0,R_1,R_2}, \mcl A = \fourpi{\emptyset}{\emptyset}{\emptyset}{R_0,0,0}, \\
&\mcl B_i = \fourpi{\emptyset}{\emptyset}{Q_2}{\emptyset}, \mcl C = \fourpi{\emptyset}{Q_1}{\emptyset}{\emptyset}, D_2 = \bmat{0\\1},
\end{align*}
where
\begin{align*}
&R_1(s,\theta) \hspace{-0.75mm}=\hspace{-0.75mm} \bmat{s-\theta&0\\0&0},~R_2(s,\theta) \hspace{-0.75mm}=\hspace{-0.75mm} \bmat{0&0\\0&\theta-s},~Q_2 =\bmat{1\\0},\\
&R_0(s) \hspace{-0.75mm}=\hspace{-0.75mm} \bmat{0&\hspace{-1.5mm}-0.1\\1&\hspace{-1.5mm}0},~Q_1(s) \hspace{-0.75mm}=\hspace{-0.75mm} \bmat{0&\hspace{-1.5mm}-\frac{s^4}{12}-\frac{s^3}{6}+\frac{s^2}{2}\\0&\hspace{-1.5mm}0}.
\end{align*}
\end{ex}
This tedious process of constructing the PIE representation has been automated in the PIETOOLS software package \cite{toolbox:pietools,manual} with a dedicated command line and GUI input formats. Typically, given a coupled ODE-PDE with sufficient boundary conditions, one can find a PIE representation of the form Eq.~\eqref{eq:PIE_full} using either the Cauchy's rule for repeated integration or PIETOOLS.

\subsection{A Side Note on $H_2$-norm Definition}
The $H_2$-norm of a system can have two different equivalent interpretations. One of the interpretations is deterministic in the sense that $H_2$-norm is defined as the system response to an impulse input (or, equivalently, an initial condition \cite{danilo}) as shown below. For this purpose, consider an abstract operator representation of a PDE with inputs $u$ and outputs $y$ given by
\begin{align}\label{eq:pde-abstract}
    \dot{\mbf x}(t) &= \mbf A\mbf x(t)+\mbf Bu(t),\quad y(t) = \mbf C \mbf x(t),
\end{align}
where $\mbf x(t)\in X\subseteq \R L_2$, $X$ is a the domain of the infinitesimal generator $\mbf A$, and $\mbf B:\R \to \R L_2$ and $\mbf C:X\to \R$ are input and output operators, respectively.
\begin{defn}\label{def:h2norm}
    Given a PDE of the form Eq.~\eqref{eq:pde-abstract}, suppose there exists $\mbf x$ differentiable for all $t\ge 0$ that satisfies Eq.~\eqref{eq:pde-abstract} for any initial condition of the form $\mbf x(0)=\mbf B u_0$ where $\mbf u_0\in \R$ and zero input (i.e., $u=0$). Then, the $H_2$-norm is given by 
    \[
    \sup\limits_{\norm{u_0}=1} \left\lbrace \norm{y}_{L_2}\;\mid \; \mat{\{\mbf x, y\}~\text{satisfy the PDE with}~u=0, \\\mbf x(0)=\mbf B u_0~, u_0\in \R }\right\rbrace.
    \]
\end{defn}

Alternatively, one can define $H_2$-norm as the largest eigenvalue of a linear operator in Hardy space, $H_2$, by using the notion of transfer functions as shown below.
\begin{defn}
    Given a PDE of the form Eq.~\eqref{eq:pde-abstract}, suppose $G$ that maps $\hat{u}$ to $\hat{y}$ given by the relation $\hat{y}(s) = G(s)\hat{u}(s)$ where $\hat y$ and $\hat u$ are Laplace transforms of the output $y$ and input $u$. Then, $H_2$-norm of the system Eq.~\eqref{eq:pde-abstract} is $\norm{G}_{H_2}$.
\end{defn}

Although both these definitions of $H_2$-norm are equivalent when $\mbf A$ generates a strongly continuous semigroup \cite{danilo}, we will use the former version (Def.~\ref{def:h2norm}) as it is more suitable for proving the results of this paper.

\section{Linear PI Inequalities}\label{sec:LPI}
We proceed by acknowledging, but without formally stating, the fact that the PDE and PIE representations are equivalent, i.e., the two representations have the same internal stability and input-output properties. See \cite{shivakumar_representation_TAC} for more details. Relying on this equivalence of input-output properties, we now formulate the $H_2$-norm optimization problems for PDEs using the PIE representation. Specifically, in this section, we present convex optimization formulations of $H_2$-norm bounding, $H_2$-optimal state estimator design, and $H_2$-optimal state feedback controller synthesis problems. 

Since PIEs are defined by PI operators, these formulations will naturally have decision variables and positivity constraints involving PI operators -- i.e., a Linear PI Inequality (LPI) problem. For example, given a PIE system $\Sigma(\mcl T, \mcl A,-,-,-)$, the following LPI is a test for stability:
\begin{align}\label{eq:stab-lpi}
\mcl P\succ 0,\quad  \mcl A^*\mcl P\mcl T+\mcl T^*\mcl P\mcl A\preceq 0.
\end{align}
Although the method to solve these problems is not described here in detail, an overview is presented below. 

In brief, these methods construct a positive PI operator using a quadratic form involving a positive matrix and \textit{$n^{th}$-order basis} of PI operators, $\mcl Z_n$. For example, $\mcl P \succeq 0$ if there exists some matrix $Q\ge0$ such that $\mcl P=\mcl Z_n^* Q\mcl Z_n=\mcl Z_n^* Q^{\half}Q^{\half} \mcl Z_n\succeq 0$, where the basis $\mcl Z_n$ is constructed using a vector of monomials in $s$ upto order $n$, $Z_n$, as
\begin{equation}\label{eq:monomial_bases}
\mcl Z_n\bmat{x\\\mbf x}(s) = \bmat{x\\Z_n(s)\mbf x(s)\\\int_a^s (Z_n(s)\otimes Z_n(\theta))\mbf x(\theta)d\theta\\\int_s^b (Z_n(s)\otimes Z_n(\theta))\mbf x(\theta)d\theta}.
\end{equation}
Thus, one can use positive matrices to parameterize positive operators $\mcl P$ and $\mcl Q$ test the feasibility of the LPI Eq.~\eqref{eq:stab-lpi} by solving the constraint $\mcl A^*\mcl P\mcl T+\mcl T^*\mcl P\mcl A = -\mcl Q$.

\subsection{LPI for $H_2$-norm Upper Bound}\label{subsec:h2norm}
Before solving the estimator and controller design problems for PIEs, we will first revisit the duality property of a PIE system. This is crucial because the duality property enables us to pose the estimation and control problems as duals of each other. Consequently, if we solve one of the two mentioned problems, we can solve the other by solving its dual. To show this duality between the two problems, we recall the following dual relation between initial conditions of $\Sigma(\mcl T,\mcl A, \mcl B,\mcl C, D)$ and its dual  $\Sigma(\mcl T^*,\mcl A^*,\mcl C^*,\mcl B^*,D^T)$ from \cite{shivakumar2020duality}. 

\begin{thm}\label{thm:intertwining}
Given $\mbf x_0$, $\bar{\mbf x}_0\in \R L_2^{m,n}$, PI operators $\mcl T, \mcl A, \mcl B, \mcl C$, and matrix $D$, if  $\{\mbf x, z\}$ satisfies $\Sigma(\mcl T,\mcl A, \mcl B,\mcl C, D)$ for initial conditions $\mcl T\mbf x_0$ and zero inputs, and $\{\bar{\mbf x},\bar z\}$ satisfies $\Sigma(\mcl T^*,\mcl A^*,\mcl C^*,\mcl B^*,D^T)$ for initial conditions $\mcl T^*\bar{\mbf x}_0$ and zero inputs, then 
\begin{align}\label{eq:intertwining}
    \ip{\mcl{T}^*\bar{\mbf x}_0}{\mbf x(t)}_{\R L_2}=\ip{\bar{\mbf x}(t)}{\mcl{T}\mbf x_0}_{\R L_2},\quad \forall ~t\ge 0.
\end{align}
\end{thm}
\begin{proof}
Proof can be found in Theorem 10 of \cite{shivakumar2020duality}.
\end{proof}

\begin{thm}\label{thm:h2-duality}
Suppose $\mcl{T}$, $\mcl{A}$, $\mcl{B}$, and $\mcl{C}$ are PI operators. Then the following statements are equivalent.
        \begin{enumerate}
                \item For any $u_0\in \R^q$, if  $\{\mbf x, z\}$ satisfies the system $\Sigma(\mcl T,\mcl A,\mcl B,\mcl C,-)$ for initial conditions $\mcl Bu_0$ and zero inputs then $\norm{z}_{L_2}\le \gamma\norm{u_0}_{\R}$.
                \item For any $\bar u_0\in \R^p$, if $\{\bar{\mbf x}, \bar{z}\}$ satisfies the system $\Sigma(\mcl T^*,\mcl A^*,\mcl C^*,\mcl B^*,-)$ for initial conditions $\mcl C^*\bar{u}_0$ and zero inputs then $\norm{\bar{z}}_{L_2}\le \gamma\norm{\bar{u}_0}_{\R}$.
        \end{enumerate}
\end{thm}
\begin{proof}
        Suppose that $\{\mbf x, z\}$ satisfy $\Sigma(\mcl T,\mcl A,\mcl B,\mcl C,-)$ for initial condition $\mcl Bu_0$ for some $u_0\in\R^q$ and zero inputs. Furthermore, let $\norm{z}_{L_2}\le \gamma\norm{u_0}_{\R}$ be valid for all $\{\mbf x, z\}$ that satisfy $\Sigma(\mcl T,\mcl A,\mcl B,\mcl C,-)$ for initial condition of the form $\mcl T\mbf x(0)=\mcl Bu_0$ for every $u_0\in \R^q$ and zero inputs.
        Let $\bar{\mbf x}(t)\in \R L_2^{m,n}$ and $\bar z(t)\in \R^q$ satisfy $\Sigma(\mcl T^*,\mcl A^*,\mcl C^*,\mcl B^*,-)$ for initial condition $\mcl C^*\bar u_0$ for some $\bar u_0\in \R^p$ and zero inputs. Then, from Eq.~\eqref{eq:intertwining} of Theorem \ref{thm:intertwining}, we have 
        \begin{align*}
                0&=\ip{\mcl{T}^*\bar{\mbf x}(0)}{\mbf x(t)}_{\R L_2}-\ip{\bar{\mbf x}(t)}{\mcl{T}\mbf x(0)}_{\R L_2}\\
                &=\ip{\mcl C^*\bar{u}_0}{\mbf x(t)}_{\R L_2}-\ip{\bar{\mbf x}(t)}{\mcl{B}u_0}_{\R L_2}\\
                &=\ip{\bar{u}_0}{\mcl C\mbf x(t)}_{\R}-\ip{\mcl B^*\bar{\mbf x}(t)}{u_0}_{\R}= \bar{u}_0^T z(t)-\bar{z}(t)^Tu_0.
        \end{align*}
        Thus, for any $t\ge 0$, we have the relationship $\bar{u}_0^T z(t) = \bar{z}(t)^Tu_0$. Squaring the left and right hand sides of the equality, we obtain
        \begin{align}\label{eq:intertwining_h2}
           (\bar{u}_0^T z(t))^T\bar{u}_0^T z(t) = (\bar{z}(t)^Tu_0)^T\bar{z}(t)^Tu_0.             
        \end{align}
        From Cauchy-Schwarz Inequality,
        \begin{align*}
            (\bar{u}_0^T z(t))^T\bar{u}_0^T z(t) \le \norm{z(t)}^2_{\R}\norm{\bar u_0}^2_{\R}.
        \end{align*}
        For any $\bar{z}(t)\in \R^p$, we know 
        \begin{align*}
                \norm{\bar{z}(t)}^2_{\R} &= \sup_{\norm{u_0}=1} \ip{\bar{z}(t)}{u_0}_{\R}^2 = \sup_{\norm{u_0}=1} (\bar{z}(t)^Tu_0)^T\bar{z}(t)^Tu_0 \\
                &\le \sup_{\norm{u_0}=1}\norm{z(t)}^2_{\R}\norm{\bar u_0}^2_{\R}.
        \end{align*}
        Hence, integrating with respect to time (from $[0,\infty)$) on both sides, we have 
        \begin{align*}
        \norm{\bar z}^2_{L_2} \le \sup_{\norm{u_0}=1}\norm{z}^2_{L_2}\norm{\bar u_0}^2_{\R} \le \gamma^2 \norm{\bar u_0}^2_{\R}.
        \end{align*}

        % Define $\bar Z_T = \int_0^T \bar{z}(t)\bar{z}(t)^T dt$ and $Z_T = \int_0^T z(t)z(t)^T dt$. Next, integrating the right hand side of Eq.~\eqref{eq:intertwining_h2} respect to time, we get
        % \begin{align*}
        %     \int_0^T (\bar{z}(t)^Tu_0)^T\bar{z}(t)^Tu_0 dt&=u_0^T \int_0^T \bar{z}(t)\bar{z}(t)^T dt u_0\\
        %     &= u_0^T \bar Z_T u_0.
        % \end{align*}
        % Again from Eq.~\eqref{eq:intertwining_h2}, we use the fact that trace of a scalar is itself and the circularity property of the trace operator, to show
        % \begin{align*}
        %         &u_0^T \bar Z_T u_0\\
        %         &= \int_0^T z(t)^T\bar u_0\bar{u}_0^T z(t)dt= \int_0^T trace(z(t)^T\bar u_0\bar{u}_0^T z(t))dt\\
        %         &= trace\left(\left(\int_0^T z(t)z(t)^T dt\right)\bar u_0\bar{u}_0^T\right)\\
        %         &\le trace(Z_T)trace(\bar u_0\bar{u}_0^T)=\norm{P_Tz}_{L_2}^2\norm{\bar{u}_0}_{\R}^2 \\
        %         &\le \gamma^2 \norm{u_0}_\R^2 \norm{\bar{u}_0}_{\R}^2,     
        % \end{align*}
        % where we also use the inequality $trace(AB)\le trace(A)trace(B)$ when $A, B$ are positive semidefinite matrices.
        % Lastly, since the terms on either side of the inequality are scalars we have
        % \begin{align*}
        % u_0^T \bar Z_T u_0 &= trace(u_0^T \bar Z_T u_0)=trace(\bar Z_T u_0u_0^T) \\
        % &\le \gamma^2 \norm{u_0}_\R^2 \norm{\bar{u}_0}_{\R}^2.
        % \end{align*}
        % Since $\bar Z_T$ is rank one, $trace(\bar Z_T)$ is the only non-zero eigenvalue of $\bar Z_T$ and hence $\norm{P_T\bar z}_{L_2}^2=trace(\bar Z_T)\le \gamma^2\norm{\bar{u}_0}^2.$ Taking limit $T\to \infty$ we have $\norm{\bar z}_{L_2}\le \gamma \norm{\bar{u}_0}_{\R}.$

        Hence, $\gamma$ is an upper bound on $H_2$-norm for $\Sigma(\mcl T^*,\mcl A^*,\mcl C^*,\mcl B^*,-)$. Since the dual and primal systems are interchangeable, necessity follows from sufficiency.
\end{proof}

As will be seen below, the above relationship leads to two formulations of $H_2$-norm analysis problem, which will later be used to solve estimator and controller design problems.

\begin{thm}\label{thm:h2-norm}
   Suppose there exist $\epsilon>0, \gamma>0$, PI operators $\mcl P\succeq \epsilon I$ and $\mcl Z$, such that either of the following two conditions hold: 
   \begin{enumerate}[label=\textbf{\thethm.\arabic*}]
   \item $\mcl T^*\mcl P\mcl A+\mcl A^*\mcl P\mcl T+ \mcl C^*\mcl C\preceq 0,\; trace(\mcl B^* \mcl P\mcl B)\le \gamma^2.$ \label{eq:h2-norm-observability}
    \item $\mcl T\mcl P\mcl A^*+\mcl A\mcl P\mcl T^*+ \mcl B\mcl B^*\preceq 0,\; trace(\mcl C \mcl P\mcl C^*)\le \gamma^2.$ \label{eq:h2-norm-controllability}
\end{enumerate}

Then, $\gamma$ is an upper bound on the $H_2$-norm of $\Sigma(\mcl T,\mcl A,\mcl B,\mcl C,-)$.
\end{thm}
\begin{proof}
    The proof of Part 1 can be found in \cite{danilo}. For Part 2, we use the equivalence of a PIE and its dual representation. Suppose there exists $\gamma, \mcl P$ that satisfy ineq.~\eqref{eq:h2-norm-controllability}. Then, from \cite[Theorem 5]{danilo}, we have that the $H_2$-norm of the PIE system $\Sigma(\mcl T^*,\mcl A^*,\mcl C^*,\mcl B^*,-)$ is upper-bounded by $\gamma$. However, from Thm.~\ref{thm:h2-duality}, we have that $H_2$-norm of a PIE and its dual are equivalent. Thus, we have that the $H_2$-norm of the PIE system $\Sigma(\mcl T,\mcl A,\mcl B,\mcl C,-)$ is upper-bounded by $\gamma$.

\end{proof}

\subsection{LPI for $H_2$-optimal Estimator}\label{subsec:h2estimator}
Given a PIE system
\begin{align}\label{eq:pie-main}
\mcl T\dot{\mbf x}(t) &= \mcl A\mbf x(t)+\mcl B_1 w(t),\notag\\
z(t) &= \mcl C_1\mbf x(t),\; y(t) = \mcl C_2 \mbf x(t)+D_{21}w(t),
\end{align}
we can design a Luenberger observer, whose dynamics is given by
\begin{align}\label{eq:pie-observer}
\mcl T\dot{\hat{\mbf x}}(t) &= \mcl A\hat{\mbf x}(t)+\mcl L(\hat y(t)-y(t)),\notag\\
\hat z(t) &= \mcl C_1\hat{\mbf x}(t),\;
\hat y(t) = \mcl C_2 \hat{\mbf x}(t)+D_{21}w(t),
\end{align}
which leads an estimation error $\mbf e = \hat{\mbf x} - \mbf x$. We can then write the dynamics of the estimation error as
\begin{align}\label{eq:error-observer}
\mcl T\dot{\mbf e}(t) &= (\mcl A+\mcl L\mcl C_2)\mbf e(t)-(\mcl B_1+\mcl LD_{21})w(t),\notag\\
\hat z(t) &= \mcl C_1\mbf e(t).
\end{align}
\begin{cor}\label{thm:h2-estimator}
Suppose there exist $\epsilon>0, \gamma>0$, matrix $W$, and PI operators $\mcl P\succeq \epsilon I$ and $\mcl Z$, such that 
\begin{align}\label{eq:h2-observe}
    &\mcl T^*\mcl P\mcl A+\mcl A^*\mcl P\mcl T+\mcl T^*\mcl Z\mcl C_2+\mcl C_2^*\mcl Z^*\mcl T + \mcl C_1^*\mcl C_1\preceq 0,\notag\\ 
    &\bmat{\mcl P&-(\mcl P\mcl B_1+\mcl ZD_{21})\\-(\mcl P\mcl B_1+\mcl ZD_{21})^*& W}\succeq 0\notag\\
    &trace(W)\le \gamma^2.
\end{align}
Then, the $H_2$-norm of $\Sigma(\mcl T,(\mcl A+\mcl L\mcl C_2),-(\mcl B_1+\mcl LD_{21}),\mcl C_1,-)$ is upper bounded by $\gamma$ where $\mcl L = \mcl P^{-1}\mcl Z$.
\end{cor}
\begin{proof}
    Suppose $\gamma, \mcl P, \mcl Z$ are as stated above. Then, from the first inequality in Eq.~\eqref{eq:h2-observe},
    \begin{align*}
        &\mcl T^*\mcl P\mcl A+\mcl A^*\mcl P\mcl T+\mcl T^*\mcl Z\mcl C_2+\mcl C_2^*\mcl Z^*\mcl T + \mcl C_1^*\mcl C_1\\
        &= \mcl T^*\mcl P(\mcl A+\mcl L\mcl C_2)+(\mcl A+\mcl L\mcl C_2)^*\mcl P\mcl T+ \mcl C_1^*\mcl C_1\le 0.
    \end{align*}
    Using Schur's Complement on the second inequality in Eq.~\eqref{eq:h2-observe}, we can say
    \[\bmat{\mcl P&-(\mcl P\mcl B_1+\mcl ZD_{21})\\-(\mcl P\mcl B_1+\mcl ZD_{21})^*& W}\succeq 0\]
    implies
    \begin{align*}
        0&\preceq W-(\mcl P\mcl B_1+\mcl ZD_{21})^*\mcl P^{-1}(\mcl P\mcl B_1+\mcl ZD_{21}) \\
        &= W-(\mcl P\mcl B_1+\mcl ZD_{21})^*\mcl P^{-1}(\mcl P\mcl B_1+\mcl ZD_{21})\\
        &= W-(\mcl B_1+\mcl LD_{21})^*\mcl P(\mcl B_1+\mcl LD_{21}),
    \end{align*}
    and $trace(W-(\mcl B_1+\mcl LD_{21})^*\mcl P(\mcl B_1+\mcl LD_{21}))\ge 0$.
    Since $trace(W)\le \gamma^2$, we have
    \[
    trace((\mcl B_1+\mcl LD_{21})^*\mcl P(\mcl B_1+\mcl LD_{21}))\le trace(W)\le \gamma^2.
    \]
    Clearly, from Part 1 of Thm.~\ref{thm:h2-norm}, we have that $\gamma$ is an upper bound on the $H_2$-norm of the PIE system $\Sigma(\mcl T,(\mcl A+\mcl L\mcl C_2),-(\mcl B_1+\mcl LD_{21}),\mcl C_1,-)$.
\end{proof}

Thus, if one solves the LPI constraints Eq.~\eqref{eq:h2-observe} (convex constraints) while minimizing the $\gamma$ (convex objective function), we can find the estimator gains $\mcl L$ whose error system has the optimal $H_2$ performance.

\subsection{LPI for $H_2$-optimal Controller}\label{subsec:h2control}
Similar to the approach taken to formulate the $H_2$-optimal estimator using Thm.~\ref{thm:h2-norm}, we can also formulate its dual problem, the $H_2$-optimal state-feedback control problem as an LPI. In this context, we consider the PIE system without observed output and distinguish between disturbance $w$ and control input $u$. Such a PIE can be represented as
\begin{align}\label{eq:pie-main-2}
\mcl T\dot{\mbf x}(t) &= \mcl A\mbf x(t)+\mcl B_1 w(t)+\mcl B_2 u(t),\notag\\
z(t) &= \mcl C_1\mbf x(t)+ D_{12}u(t).
\end{align}
The goal is to find a state-feedback $u(t) = \mcl K\mbf x(t)$, such that the $H_2$ performance of the closed-loop system given by $\Sigma(\mcl T,(\mcl A+\mcl B_2\mcl K),\mcl B_1,(\mcl C_1+\mcl D_{12}\mcl K),-)$ is optimal. Since controller synthesis is the dual problem for $H_2$-optimal estimator problem, we use the dual version of LPI in Thm.~\ref{thm:h2-norm} (i.e., Eq.~\eqref{eq:h2-norm-controllability}) to formulate the LPIs for controller synthesis as shown below.

\begin{cor}\label{thm:h2-controller}
Suppose there exist $\epsilon>0, \gamma>0$, matrix $W$, and PI operators $\mcl P\succeq \epsilon I$ and $\mcl Z$, such that 
\begin{align}\label{eq:h2-control}
    &\mcl A\mcl P\mcl T^*+\mcl T\mcl P\mcl A^*+\mcl B_2\mcl Z\mcl T^*+\mcl T\mcl Z^*\mcl B_2^* + \mcl B_1\mcl B_1^*\preceq 0,\notag\\ 
    &\bmat{\mcl P&(\mcl C_1\mcl P+D_{12}\mcl Z)^*\\\mcl C_1\mcl P+D_{12}\mcl Z& W}\succeq 0\notag\\
    &trace(W)\le \gamma^2.
\end{align}
Then, the $H_2$-norm of $\Sigma(\mcl T,(\mcl A+\mcl B_2\mcl K),\mcl B_1,(\mcl C_1+\mcl D_{12}\mcl K),-)$ is upper bounded by $\gamma$ where $\mcl K = \mcl Z\mcl P^{-1}$.
\end{cor}
\begin{proof}
    The proof follows an approach similar to the proof of Thm.~\ref{thm:h2-estimator}. Hence, we will only provide an outline here. Suppose $\gamma, \mcl P,\mcl Z$ are as stated above. Then, from the first inequality in Eq.~\eqref{eq:h2-control}, we have
    \begin{align*}
        \mcl T\mcl P(\mcl A+\mcl B_2\mcl K)^*+(\mcl A+\mcl B_2\mcl K)\mcl P\mcl T^*+ \mcl B_1\mcl B_1^*\le 0.
    \end{align*}
    From second inequality in Eq.~\eqref{eq:h2-control},
    \begin{align*}
        W-(\mcl C_1+D_{12}\mcl K)\mcl P(\mcl C_1+D_{12}\mcl K)^*\succeq 0.
    \end{align*}
    Hence, $trace((\mcl C_1+D_{12}\mcl K)\mcl P(\mcl C_1+D_{12}\mcl K)^*)\le \gamma^2$. Thus, the $H_2$-norm of the PIE system $\Sigma(\mcl T,(\mcl A+\mcl B_2\mcl K),\mcl B_1,(\mcl C_1+\mcl D_{12}\mcl K),-)$ is upper bounded by $\gamma$.
\end{proof}

Again, similar to solving for $H_2$-optimal estimator, we can solve the LPI constraints Eq.~\eqref{eq:h2-control} while minimizing $\gamma$ to find the controller gains $\mcl K$ and obtain the closed-loop system that has the optimal $H_2$ performance. However, to find the gains (both estimator and controller) requires inversion of a PI operator $\mcl P$ --- an iterative approach was described in \cite{shivakumar2020duality} and will not be discussed here.

% \subsection{Inversion of Positive PI operator}

\section{Numerical Examples}\label{sec:numerical}
We will use the Matlab toolbox that was developed to solve LPI optimization problems, PIETOOLS, because the toolbox offers convenient Matlab functions to convert PDEs to PIE, declare PI decision variables, add LPI constraints, and solve the resulting optimization problem. We refer to the PIETOOLS User Manual~\cite{manual} for details. For the following two examples, namely, an unstable reaction-diffusion PDE and a neutrally-stable Euler-Bernoulli beam, we will use PIETOOLS toolbox to obtain PIE representation of the PDEs. Then, we will apply the results from Cor.~\ref{thm:h2-estimator} and \ref{thm:h2-controller} to find the $H_2$-optimal observers and controllers. Utilizing the helper functions in PIETOOLS to invert positive PI operators and we construct the closed-loop observer and controller systems, which are simulated using first-order backward difference integration scheme for certain initial conditions and no disturbance. For each example, we also provide a \textit{numerical estimate} of the $H_2$-norm (i.e., $\norm{z}_{L_2}/\norm{u_0}_{\R}$) as observed in the simulations --- i.e., by performing numerical integration of the simultation output $z(t)^2$ to obtain $\norm{z}_{L_2}$.

\subsection{Estimation and Control of Reaction-diffusion PDE}\label{subsec:react-diff}
\begin{figure}[!t]
    \centering
    \includegraphics[width=0.5\textwidth]{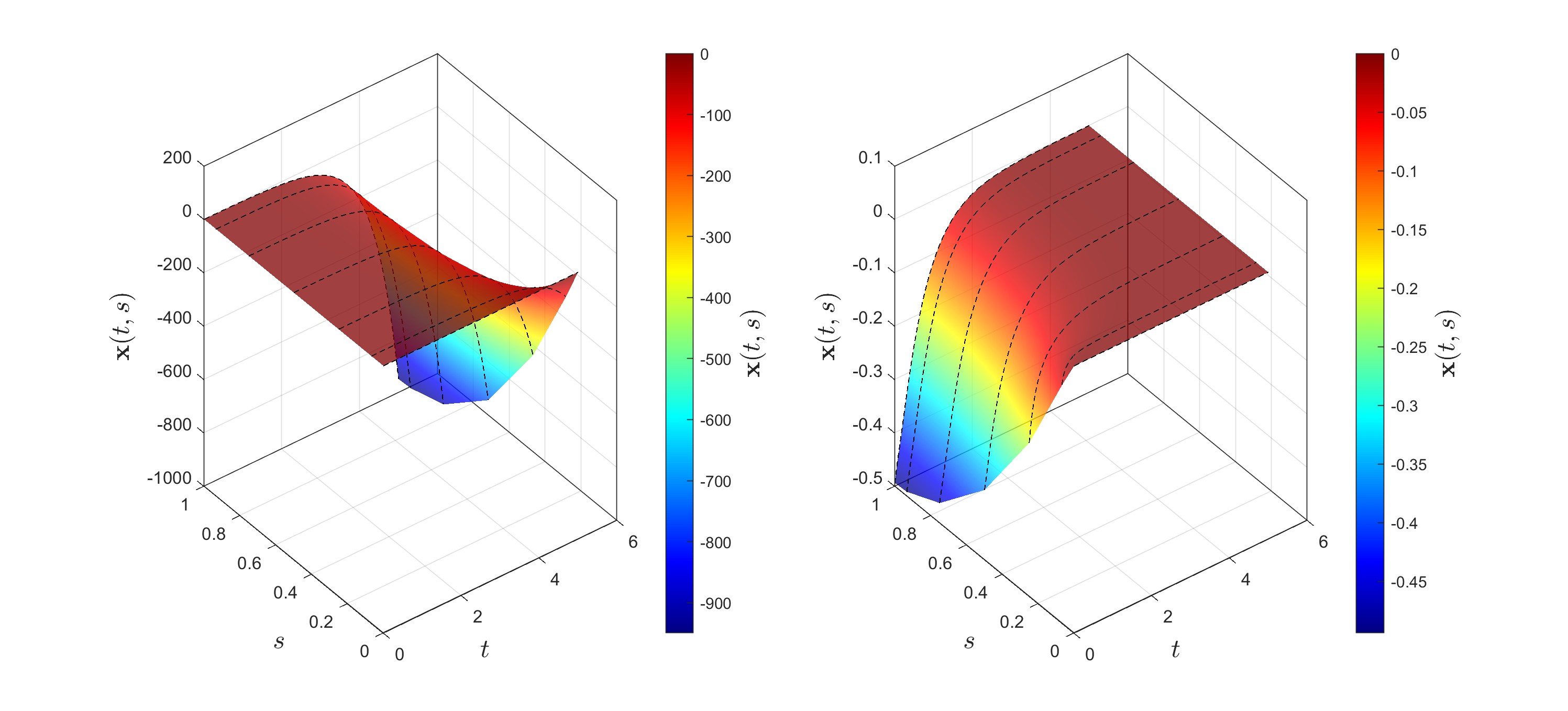}
    \caption{This figure plots the response of the system Eq.~\eqref{eq:react-diff} without control (on the left) and with control input (on the right) under zero disturbance and initial conditions $\mbf x(0,s)=\frac{s^2-2s}{2}$ while considering the regulated output as defined in Eq.~\eqref{eq:react-diff}. $H_2$-norm of the closed-loop system is $1.79$ (\textit{numerical estimate} $1.21$).}
    \label{fig:control-react-diff}
\end{figure}

In this example, we consider the reaction-diffusion PDE given by
{\small
\begin{align}\label{eq:react-diff}
    &\dot{\mbf x}(t,s) = 3\mbf x(t,s)+(s^2+0.2)\partial_s^2\mbf x(t,s)+\frac{s^2-2s}{2}w(t)+u(t),\notag\\
    &z(t)= \bmat{\int_0^1 \mbf x(t,s)ds\\u(t)}, \quad y(t) = \mbf x(t,1)+w(t),\notag\\ 
    &\mbf x(t,0)=\partial_s \mbf x(t,1)=0.
\end{align}}
For the above PDE, we use PIETOOLS toolbox to obtain a PIE representation and then solve the LPI optimization problems in Cor.~\ref{thm:h2-estimator} and \ref{thm:h2-controller} to obtain the gains corresponding to the $H_2$-optimal estimator and state-feedback controller. Then, the closed-loop PIE system is constructed and simulated using the PIESIM module of the PIETOOLS toolbox in MATLAB to find the system's response under a zero disturbance and initial conditions $\mbf x(0,s)=\frac{s^2-2s}{2}$ (i.e., $u_0=1$). The initial-condition response of the PDE.~\eqref{eq:react-diff} (without control and with $H_2$-optimal state-feedback control) are presented in Fig.~\ref{fig:control-react-diff}. The $H_2$-norm bound obtained by solving the LPI in Thm.~\ref{thm:h2-controller} is stated in the caption along with the \textit{numerical estimate}.

\begin{figure}[!t]
    \centering
    \includegraphics[width=0.5\textwidth]{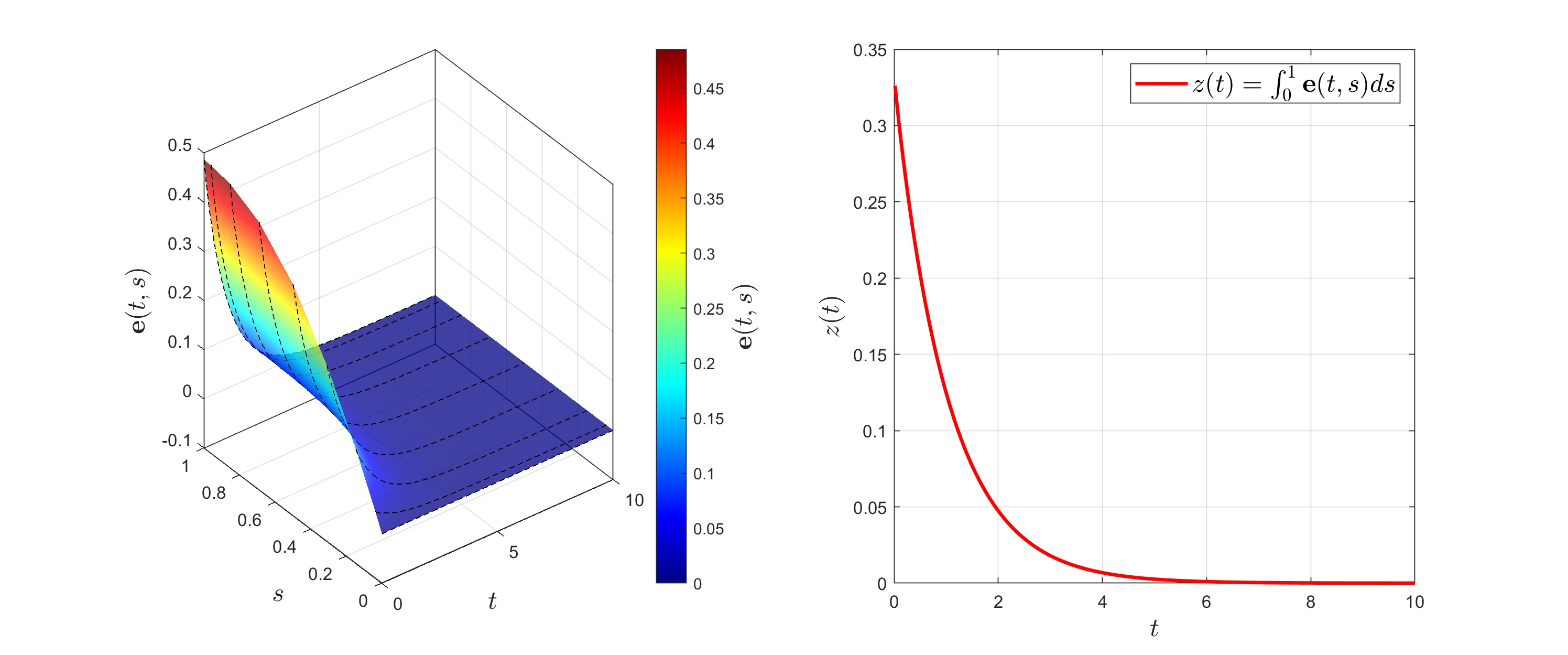}
    \caption{This figure plots the error ($\mbf e=\hat{\mbf x}-\mbf x$) in the state estimate (on the left) and regulated error output (on the right) of the state observer for the system Eq.~\eqref{eq:react-diff}. The observer is initialized with zero initial conditions, whereas the PDE state starts with an initial condition $\mbf x(0,s)=\frac{s^2-2s}{2}$ and is under zero disturbance. $H_2$-norm of the observer error system is $1.37$ (\textit{numerical estimate} $0.23$) where regulated and observed outputs are as defined in Subsec.~\ref{subsec:react-diff}.}
    \label{fig:observer-react-diff}
\end{figure}

For the observer simulation, we initialize the observer state at zero, while the PDE state is initialized as earlier. In Fig.~\ref{fig:observer-react-diff}, we only show the response of the error system --- i.e., we plot the error between state-estimate ($\hat{\mbf x}$) and actual state ($\mbf x$), given by $\mbf e = \hat{\mbf x}-\mbf x$. Additionally, we also plot the regulated output of the error system given by $\hat{z}(t)-z(t) = \int_0^1 \mbf e(t,s)ds.$

\subsection{Estimation and Control of Euler-Bernoulli Beam Equation}\label{subsec:eb}
\begin{figure}[!t]
    \centering
    \includegraphics[width=0.5\textwidth]{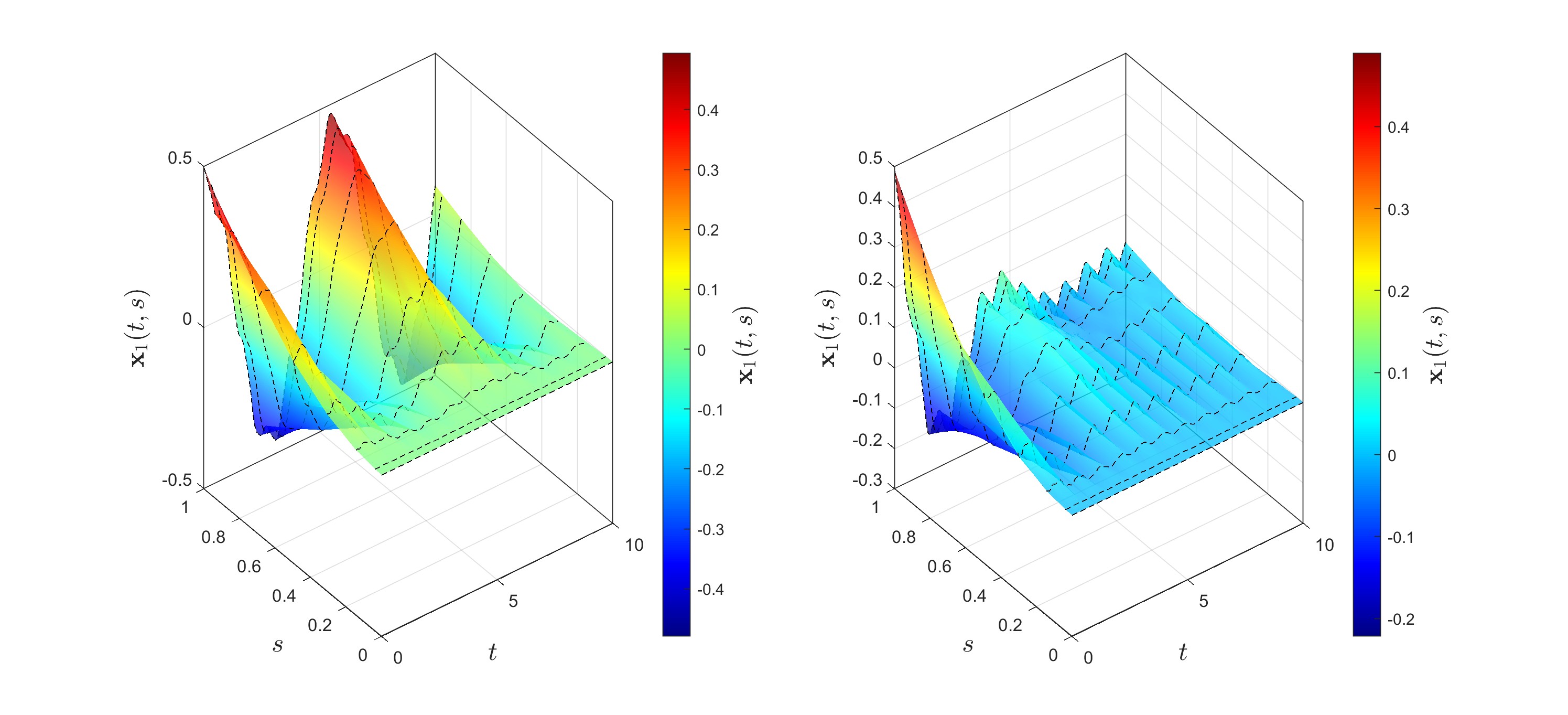}
    \caption{This figure plots the response of the system Eq.~\eqref{eq:eb} without control (on the left) and with control input (on the right) under zero disturbance and initial conditions $\mbf x(0,s)=col(-0.5s^2,0)$ while considering the regulated output as defined in Eq.~\eqref{eq:eb}. $H_2$-norm of the closed-loop system is $0.78$ (\textit{numerical estimate} $0.29$).}
    \label{fig:control-eb}
\end{figure}

In this example, we consider the Euler-Bernoulli beam equation introduced in Ex.~\ref{ex:EB-representation},
\begin{align}\label{eq:eb}
&\dot{\mbf{x}}(t,s) = \bmat{0&-0.1\\1&0}\partial_s^2\mbf{x}(t,s)+\bmat{-0.5s^2\\0}w(t)+\bmat{1\\0}u(t),\notag\\
&\bmat{1&0}\mbf x(t,0) = \bmat{1&0}\partial_s\mbf x(t,0) = 0,\notag\\
&\bmat{0&1}\mbf x(t,1) =\bmat{0&1}\partial_s\mbf x(t,1) = 0,\notag\\
&z(t) = \bmat{\int_0^1 \bmat{0&1-s-s^2}\mbf x(t,s)ds\\u(t)},\notag\\ 
&y(t) = \bmat{1&0}\mbf x(t,1)+w(t).
\end{align}
We obtain PIE representation as stated earlier and then solve the LPI optimization problems in Cor.~\ref{thm:h2-estimator} and \ref{thm:h2-controller} to obtain the gains corresponding to the $H_2$-optimal estimator and state-feedback controller. Then, the PIEs are simulated using the PIESIM to find the system's response under zero disturbance $w$ and initial conditions $\mbf x(0,s) = col(-0.5s^2,0)$ (i.e., $u_0=1$). Similar to the reaction-diffusion PDE example, we plot initial-condition response of the PDE without control and with $H_2$-optimal state-feedback control in Fig.~\ref{fig:control-eb}. The $H_2$-norm bound from Thm.~\ref{thm:h2-controller} and the \textit{numerical estimate} from the simulation are stated in the caption.

\begin{figure}[!t]
    \centering
    \includegraphics[width=0.5\textwidth]{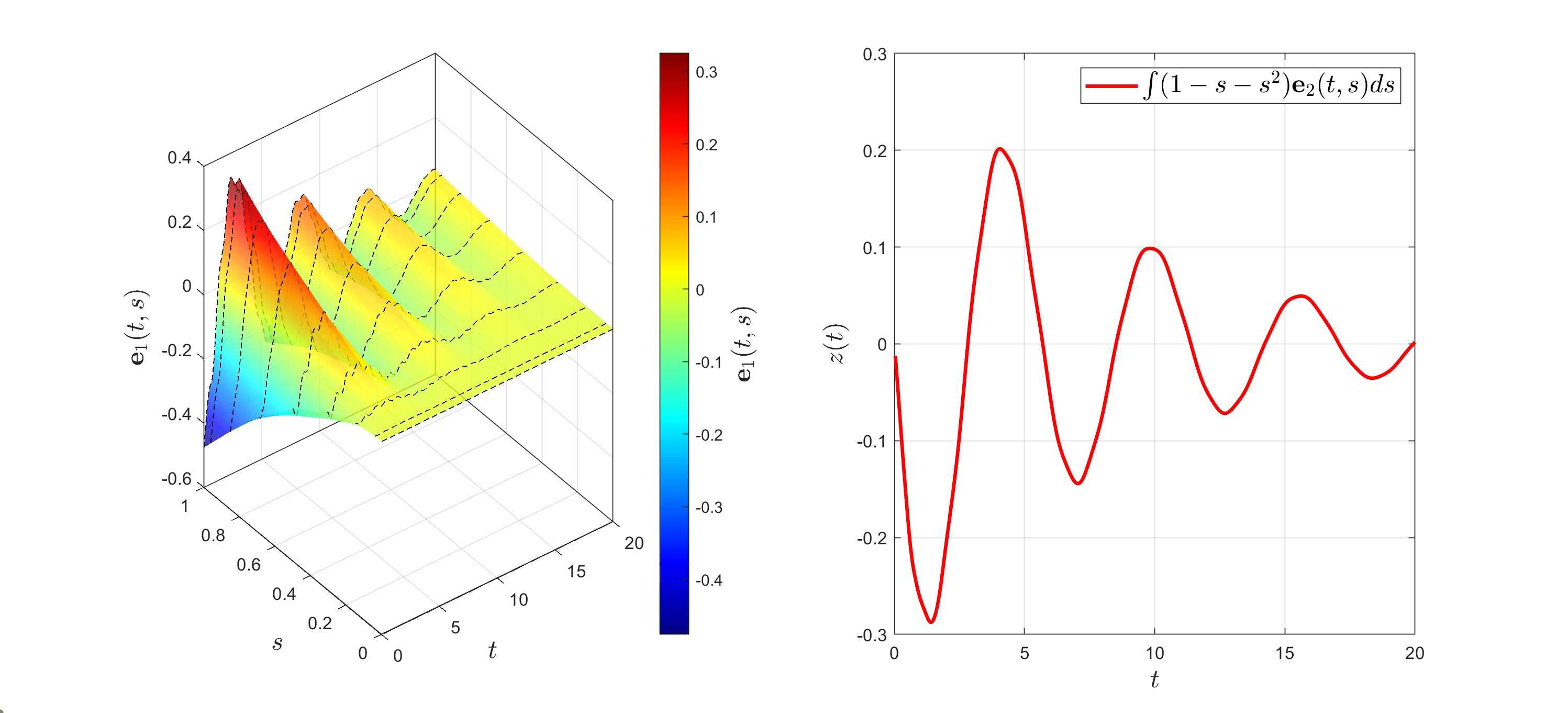}
    \caption{This figure plots the first component of the error ($\mbf e_1=\hat{\mbf x}_1-\mbf x_1$) in the state estimate (on the left) and regulated error output (on the right) of the state observer for the system Eq.~\eqref{eq:eb}. The observer is initialized with zero initial conditions, whereas the PDE state starts with an initial condition $\mbf x(0,s)=col(-0.5s^2,0)$ and is under zero disturbance. $H_2$-norm of the observer error system is $0.57$ (\textit{numerical estimate} $0.48$) where regulated and observed outputs are as defined in Subsec.~\ref{subsec:eb}.}
    \label{fig:observer-eb}
\end{figure}

In Fig.~\ref{fig:observer-eb}, we see the response of the error system for PDE initial conditions stated above and zero observer initial condition. In this case, we plot the error in only one of the states and the regulated output in the other, given by $\hat{z}(t)-z(t) = \int_0^1 (1-s-s^2)\mbf e_2(t,s)ds$.

\section{Conclusion}\label{sec:conclusion}
In this paper, we solved the $H_2$-optimal estimation and control problems for PDEs using the PIE framework developed for the analysis and control of PDE systems. Since formulating a PDE analysis/control problem using the PIE representation does not introduce any conservatism and leads to solvable convex optimization problems called Linear PI Inequalities (LPIs), we showed that $H_2$ analysis, estimation and control problems for PDEs can be solved using convex optimization without conservatism. For this purpose, we utilized an alternative definition of $H_2$-norm of a system that does not rely on the transfer function or impulse input; Instead, we characterized $H_2$-norm as the gain from an initial condition to the output of the system. Using this alternative, but equivalent, definition of $H_2$-norm, we showed that a PIE system and its corresponding dual PIE system have the same $H_2$-norm. Using this duality, we formulated two versions of LPIs problems (a primal and dual) to upper bound the $H_2$-norm of the system that were later used to formulate $H_2$-optimal state estimator and state-feedback control problems for PIEs as convex LPI optimization problems. By solving these LPI optimization problems, we demonstrated the application of this framework in estimator design and controller synthesis for PDE numerical examples. 

\addtolength{\textheight}{-7cm}

\bibliography{references}
\bibliographystyle{IEEEtran}
\end{document}